\pdfoutput=1% forces arxiv to use pdflatex to process pdf plots & figures
\documentclass[a4paper,12pt]{article}
\usepackage{ifpdf}
\ifpdf
\usepackage[pdftex]{graphicx}
\usepackage[usenames,dvipsnames]{color}
\else
\usepackage{graphicx}
\usepackage[usenames,dvips]{color}
\fi
\usepackage{amsmath}
\usepackage{amssymb}
\usepackage{amsthm}
\usepackage[colorlinks=true]{hyperref}   % for electronic versions
\usepackage{caption}
\usepackage{paper}

\numberwithin{equation}{section}

\begin{document}

\title{Two dimensional axisymmetric\\ smooth lattice Ricci flow.}
\author{%
Leo Brewin\\[10pt]%
School of Mathematical Sciences\\%
Monash University, 3800\\%
Australia}
% \date{21-Oct-2015}% 1st version started
% \date{12-Nov-2015}% 1st version finished
\date{26-Nov-2015}% 2nd version finished
% \date{29-Jan-2015}% finished
\reference{Preprint}
% \reference{Preprint: arXiv:1505.00067\\[5pt]
%            Journal: {\it Class.Q.Grav}  {\bf 32} (2015) 195008}

\maketitle

\begin{abstract}
\noindent
A lattice based method will be presented for numerical investigations of Ricci flow. The
method will be applied to the particular case of 2-dimensional axially symmetric initial
data on manifolds with $S^2$ topology. Results will be presented that show that the method
works well and agrees with results obtained using contemporary finite difference methods.
\end{abstract}

% ============================================================================================
\section{Introduction}
\label{sec:intro}

The Ricci flow \cite{hamilton:1982-01} of a metric $g$ is described by the equation
\begin{equation}
   \label{eqn:RicciFlow}
   \frac{\partial g}{\partial t} = -2Ric(g)
\end{equation}
where $Ric(g)$ is the Ricci tensor of the metric $g$. Though there is an extensive
literature on the mathematical properties of Ricci flows
\cite{chow:2004-01,chow:2006-01,morgan:2007-01,topping:2006-01} there is far less material
concerning numerical methods for Ricci flow
\cite{rubinstein:2005-01,garfinkle:2005-01,miller:2014-01}.

In 2005 Rubinstein and Sinclair \cite{rubinstein:2005-01} presented a some numerical studies
of 2-dimensional Ricci flow for two classes of geometries, one in which the geometry was
axisymmetric, and a second more general class in which the axisymmetric condition was
dropped. In both cases the surfaces were closed 2-surfaces topologically equivalent to a
2-sphere. They displayed their results for a variety of initial configurations and where
possible, displayed the evolved 2-surfaces as isometric embeddings in $E^3$. It is well
known that not all 2-geometries can be realised as an isometric embedding in $E^3$. However,
for the axisymmetric case they showed that the Ricci flow would preserve the existence of an
embedding in $E^3$. They also showed, for the non-axisymmetric case, that under the Ricci
flow it may not be possible, at later times, to embedded the surface isometricly in $E^3$.

For both classes of geometries Rubinstein and Sinclair were able to demonstrate the expected
behaviour for late term evolution under Ricci flow, namely, that in the absence of
neck-pinching, the geometry evolves towards that of a 2-sphere. Their numerical studies were
however limited by numerical issues that caused the evolution to fail at late times.

Other works include a series of papers by Kim and co-workers \cite{jin:2008-01,jin:2007-01}
in which they used a modified Ricci flow to solve the problem of finding a 2-metric that
matches a prescribed Gaussian curvature. Miller and his co-workers \cite{miller:2014-01}
have used their extensive experience in the Regge calculus to explore numerical simulations
of Ricci flow in $S^3$. They have also recently reported on neck pinching singularities in
3-dimensions \cite{alsing:2013-01}. The works of Kim and Miller both employ simplicial
methods to estimate the curvatures. In contrast Garfinkle and Isenberg
\cite{garfinkle:2005-01} used traditional finite difference methods to follow the Ricci flow
for a class of geometries on $S^3$. They were particularly interested in the formation of
neck pinching singularities and the development of critical behaviour in the solutions along
the lines found by Choptuik \cite{choptuik:1993-01}.

The purpose of this paper is to use the axisymmetric models of Rubinstein and Sinclair as a
basis to demonstrate a novel method for computational differential geometry. This method is
based on the idea that a smooth geometry can be approximated by a finite set of overlapping
cells, with each cell defined by a small set of vertices and legs, and in which each cell
carries a locally smooth metric. Such a structure is known as a smooth lattice and has been
successfully employed in constructing numerical solutions of the Einstein field equations
(see \cite{brewin:2014-01,brewin:2010-02,brewin:2010-03}).

Since this smooth lattice method is not widely known it seems appropriate to provide a short
review of the method.
%That review is provided in the following section.

\testbreak

% ============================================================================================
\section{A smooth lattice}
\label{sec:smoothlattice}

A smooth lattice can be considered as a generalisation of a simplicial lattice this being a
manifold built from a finite collection of simplices each endowed with a flat metric. The
simplicial lattice, also known as a piecewise flat lattice, is the basis of many studies in
discrete differential geometry and has also been used in General Relativity in a formulation
known as the Regge Calculus (\cite{regge:1961-01,gentle:2002-01,williams:2009-01}). The
simplicial and smooth lattices differ in one important aspect -- the smooth lattice uses a
locally smooth metric rather than a piecewise flat metric. The domain over which each
locally smooth metric is defined is bounded by the nearby vertices and legs. Each such
domain is known as a computational cell and pairs of adjacent cells are allowed to have a
non-trivial overlap.

The picture that emerges form this is much like that of an atlas of charts that covers a
manifold. Each chart can be viewed as a computational cell while the overlap between pairs
of cells defines the transition functions between overlapping charts.

\begin{figure}[ht]
\Figure{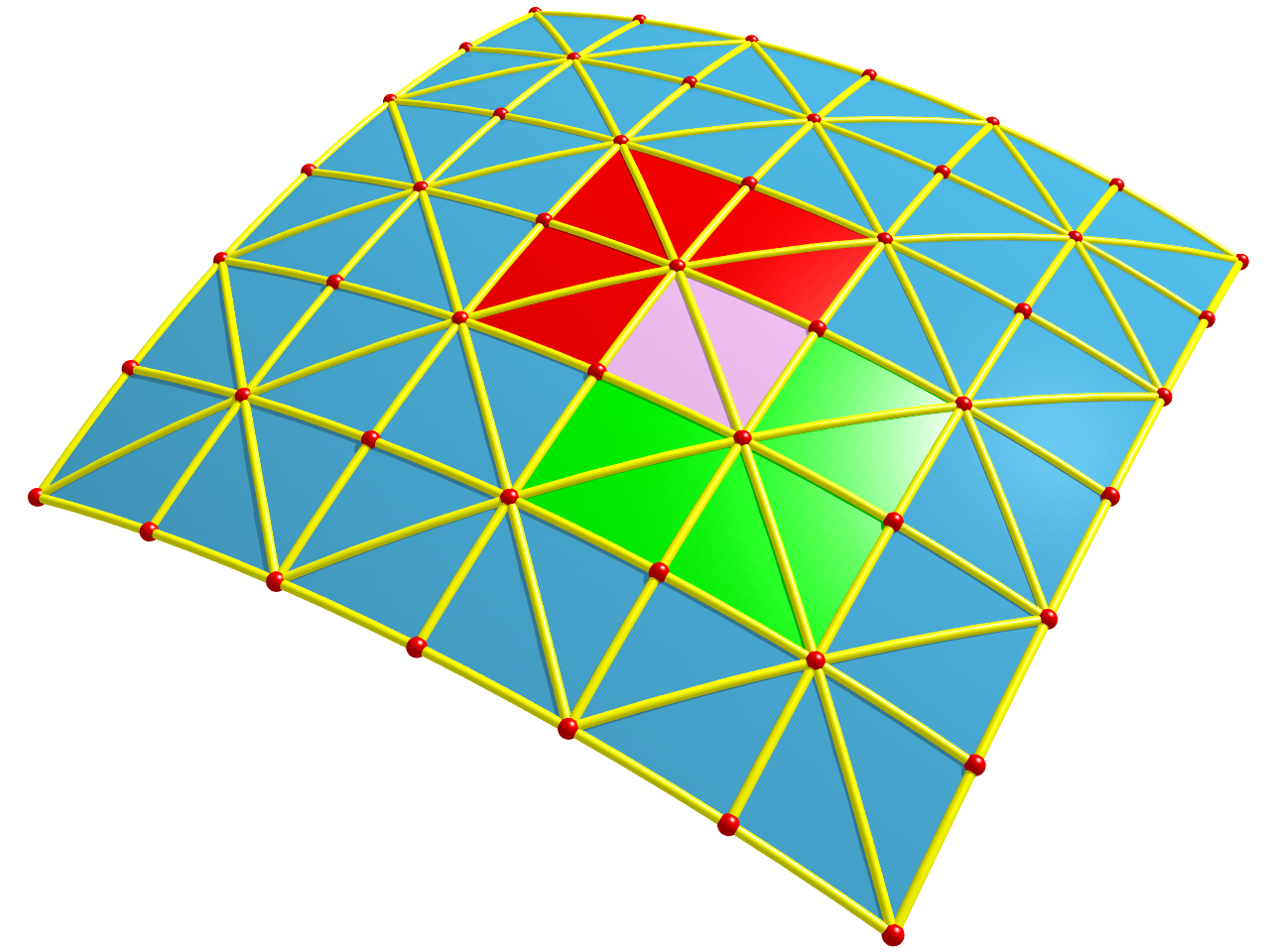}{0.90}
\vskip 15pt
\caption{\normalfont%
An example of a 2-dimensional lattice. This shows two overlapping computational cells, one
red, the other green. The two cells share data in the pink region. Note that cells in this
diagram were chosen to form a regular structure for purely aesthetic reasons. On a general
lattice the number of triangles meeting at a vertex may vary from vertex to vertex and
likewise for the leg-lengths. Note that in the smooth lattice method each leg in the lattice
is taken to be a short geodesic segment of the smooth geometry for which the lattice is an
approximation.}
\label{fig:2dlattice}
\end{figure}

A simple example of 2-dimensional smooth lattice is shown in figure (\ref{fig:2dlattice}).
This shows a lattice composed of the yellow legs and red vertices that forms a discrete
approximation to the underlying smooth blue surface. The figure also shows two overlapping
computational cells.

The data that is employed in a smooth lattice is at least the set of leg lengths and
possibly other geometric data such as the Riemann tensor within each computational cell.
Note that the leg length assigned to any one leg is a property of that leg (i.e., if a leg
is shared by a pair of cells then each cell uses the same length for that leg). It should
also be noted that each leg is viewed as a geodesic segment of a (possibly unknown) smooth
geometry (the geometry for which the smooth lattice is an approximation). Without this
assumption the path joining the vertices of a leg would be ambiguous (and cell dependent).
This geodesic assumption imposes a soft constraint on the vertices of the lattice -- they
must be chosen so that the geodesic for each leg is uniquely defined. For smooth geometries
this is always possibly by a suitable refinement of the lattice.

Using a local smooth metric allows the usual machinery of differential geometry to be
applied directly to the lattice. In contrast, the piecewise constant nature of the metric on
a simplicial lattice requires considerable care when applying differential operators to the
lattice. For example, the curvature tensor on a simplicial lattice must be interpreted in the
sense of distributions.

Since each leg in a computational cell is assumed to be defined as the unique geodesic for
that leg, it follows that it is always possible to construct a local set of Riemann normal
coordinates in each computational cell. Denote these coordinates by $x^\mu$. Then the smooth
metric, in this cell, can be written as
\begin{equation}
   \label{eqn:rncMetric}
   g_{\mu\nu}(x) = g_{\mu\nu}
                 + \frac {1}{3} R_{\mu\alpha\nu\beta} x^\alpha x^\beta +\BigO{L^3}
\end{equation}
where $g_{\mu\nu} = \diag(1,1,1,\cdots)$ and $L$ is a typical length scale for the cell.

Using Riemann normal coordinates allows some useful quantities to be easily computed. For
example the arc-length $L_{ij}$ of the geodesic connecting vertices $i$ and $j$ is given by
\begin{equation}
   \label{eqn:rncLsq}
   L^2_{ij} = g_{\mu\nu} \Delta x^\mu_{ij} \Delta x^\nu_{ij}
            - \frac {1}{3} R_{\mu\alpha\nu\beta} x^\mu_i x^\nu_i x^\alpha_j x^\beta_j
            + \BigO{L^5}
\end{equation}
while the unit tangent vector $v^\mu$ to the geodesic, at vertex $i$, is given by
\begin{equation}
   \label{eqn:rncGeodesicTangent}
  v^\mu
  = \frac{1}{L_{ij}}
    \left(\Delta x^\mu_{ij}
          - \frac{1}{3}R^\mu{}_{\alpha\beta\rho}
             x^\rho_i\Delta x^\alpha_{ij}\Delta x^\beta_{ij}\right)
          + \BigO{L^4}
\end{equation}
In the above pair of equations $x^\mu_i$ are the Riemann normal coordinates of vertex $i$ and
$\Delta x^\mu_{ij} = x^\mu_i - x^\mu_j$.

Both of these equations will be used later when developing one particular method to evolve
the Ricci flow.

% LCB: Comment on the common practice to use normal coordinates for local calculations
% and how the lattice stitches these together to get a global result. Or something on that line.
% See my paper on Regge-vs-SLGR.

% ============================================================================================
\section{Mathematical formulation}
\label{sec:theEquations}

% --------------------------------------------------------------------------------------------
\subsection{Rubinstein and Sinclair}
\label{sec:RubinSinc}

For their axisymmetric geometries Rubinstein and Sinclair choose coordinates $(\rho,\theta)$
in which the 2-metric is given by
\begin{equation}
   \label{eqn:RS2metric}
   ds^2 = h(\rho) d\rho^2 + m(\rho)d\theta^2
\end{equation}
with $0\le\rho\le\pi$ and $0\le\theta<2\pi$. Smoothness of the metric at $\rho=0$ and
$\rho=\pi$ requires
\begin{equation}
   \label{eqn:RSsmooth}
   \left(\sqrt{m(\rho)}\right)' = \sqrt{h(\rho)}\qquad\text{at }\rho=0,\pi
\end{equation}
where the dash $'$ denotes the derivative with respect to $\rho$.

The components of the Ricci tensor are, for $0<\rho<\pi$,
\begin{equation}
   \begin{aligned}
      R_{\theta\rho} &= R_{\rho\theta} = 0\\[5pt]
      R_{\rho\rho} & = \frac{(m')^2}{4m^2} - \frac{m''}{2m} + \frac{m'h'}{4mh}\\[5pt]
      R_{\theta\theta} &= \frac{(m')^2}{4mh} - \frac{m''}{2h} + \frac{m'h'}{4h^2}
   \end{aligned}
\end{equation}
while at the the poles, where $\rho=0$ or $\rho=\pi$,
\begin{equation}
   \begin{aligned}
      R_{\theta\theta} &= R_{\theta\rho} = R_{\rho\theta} = 0\\[5pt]
      R_{\rho\rho} &= \frac{m''''}{4m''} - \frac{h''}{2h}
   \end{aligned}
\end{equation}

The Ricci flow equation (\ref{eqn:RicciFlow}), for this metric in these coordinates, can
then be written as
\begin{equation}
   \label{eqn:RSFlowA}
   \begin{aligned}
      \frac{\partial m}{\partial t} &= -\frac{(m')^2}{2mh}
                                       + \frac{m''}{h} - \frac{m'h'}{2h^2}\\[5pt]
      \frac{\partial h}{\partial t} &= -\frac{(m')^2}{2m^2}
                                       + \frac{m''}{m} - \frac{m'h'}{2mh}
   \end{aligned}
\end{equation}
for $0<\rho<\pi$ and as
\begin{equation}
   \label{eqn:RSFlowB}
   \begin{aligned}
      \frac{\partial m}{\partial t} &= 0 \\[5pt]
      \frac{\partial h}{\partial t} &= \frac{h''}{h} - \frac{m''''}{2m''}
   \end{aligned}
\end{equation}
at the poles.

% --------------------------------------------------------------------------------------------
\subsection{Smooth lattice Ricci flow}
\label{sec:SLRF}

The basic data on a smooth lattice includes the leg-lengths, the Riemann normal coordinates
and the corresponding Riemann curvatures all of which must be considered as functions of
time under the Ricci flow. The question then must be -- how might each of these quantities
be evolved forward in time?

Consider first the evolution of the leg-lengths. The natural starting point would be the
basic Ricci flow equation which, in coordinate form, can be written as
\begin{equation}
   \label{eqn:RicciFlowCoord}
   \frac{\partial g_{\mu\nu}}{\partial t} = -2 R_{\mu\nu}
\end{equation}
Consider a typical leg defined by the vertices $i$ and $j$. Sub-divide this leg into many
short segments and let a typical segment have end points $a$ and $b$. Our first step will be
to form an estimate for the evolution of the segment followed by a summation over all
segments to obtain an evolution equation for the leg.

Choose a new set of coordinates $y^\mu$ tied to the vertices of the segment (i.e., the
coordinates of the vertices do not change with time). The squared length of the segment
joining $a$ and $b$ can be estimated as $L^2_{ab} = g_{\mu\nu}\Delta y^\mu_{ab}\Delta
y^\nu_{ab} + f_{ab}L^3_{ab}$ where $\Delta y^\mu_{ab} = y^\mu_b - y^\mu_a$ and
$f_{ab}L^3_{ab}$ is the truncation error (in using the first term to estimate $L^2_{ab}$).
Note that if the geodesic joining $i$ to $j$ does not pass through any curvature
singularites then $f_{ab}$ will be finite along that geodesic. Let $v^\mu_{ab}$ be a unit
vector tangent to the geodesic pointing from $a$ to $b$. Then from
(\ref{eqn:RicciFlowCoord}) and using $\Delta y^\mu_{ab} = v^\mu_{ab} L_{ab}$, it follows
that
\begin{equation}
   \frac{\partial L_{ab}}{\partial t} = - R_{\mu\nu} v^\mu_{ab} v^\nu_{ab} L_{ab}
                                      + \frac{1}{2} f_{ab}L^2_{ab}
\end{equation}
which, upon summing over all segments, leads to
\begin{equation}
   \frac{\partial L_{ij}}{\partial t}
   = - \sum_{ab} R_{\mu\nu} v^\mu_{ab} v^\nu_{ab} L_{ab}
     + \frac{1}{2}\sum_{ab} f_{ab} L^2_{ab}
\end{equation}
The final step is to take the limit as the number of segments $N$ approaches infinity. The
first sum on the right hand side is in the form of a Riemann sum while the second term is
subject to
\begin{equation}
   0 \le \Big\vert \sum_{ab} f_{ab} L^2_{ab}\Big\vert
     \le N\max_{ab}\vert f_{ab}\vert L^2_{ab}
\end{equation}
and thus vanishes as $N\rightarrow\infty$ provided $\max_{ab} L_{ab} = \BigO{1/N}$. This
leads to the basic evolution equation for the leg-lengths, namely
\begin{equation}
   \label{eqn:evolveLij}
   \frac{\partial L_{ij}}{\partial t} = - \int_i^j R_{\mu\nu} v^\mu v^\nu \> ds
\end{equation}
where $s$ is the arc-length along the geodesic joining $i$ to $j$ and $v^\mu(s)$ is the
unit-tangent to the geodesic. Note that since the integrand is a pure scalar the
restriction to the adapated coordinates $y^\mu$ is no longer required. Thus equation
(\ref{eqn:evolveLij}) is valid in all frames.

In any numerical code the integral in (\ref{eqn:evolveLij}) must of course be estimated in
terms of data available on the lattice. One example would be a Trapezoidal rule such as
\begin{equation}
   \label{eqn:Trapezoidal}
   \int_i^j R_{\mu\nu} v^\mu v^\nu \> ds
   \approx \frac{1}{2}
     \left( \left(R_{\mu\nu} v^\mu v^\nu\right)_i
   +        \left(R_{\mu\nu} v^\mu v^\nu\right)_j \right) L_{ij}
\end{equation}
in which the values for $R_{\mu\nu}$ are obtained from two cells, one based on vertex $i$,
the other based on vertex $j$. All of the results presented in section (\ref{sec:results})
used the above approximation. Note also that for the case of a 2-metric,
$2R_{\mu\nu} = R g_{\mu\nu}$, and thus the right hand side of (\ref{eqn:evolveLij}) can be
expressed as $-\int_i^j Rds/2$.

This shows clearly that if the Riemann curvatures are known then the leg-lengths can be
evolved. So the question now is -- how are the Riemann curvatures obtained? Two methods will
be presented, one in which the Riemann curvatures are derived from the leg-lengths and a
second method in which they are evolved in consort with the leg-lengths. In both cases the
Riemann normal coordinates are obtained directly from the lattice data.

The first method, in which the curvatures are derived from the leg-lengths, works as
follows. Consider a typical 2-dimensional cell consisting of a set of $N$
triangles\footnote{Note that the use of triangles is simply for the sake of argument, a
2-dimensional smooth lattice could also be built from any other shape provided that the
computations for the curvatures and coordinates are well defined.} sharing a common vertex.
The data to be computed are the $2N+2$ coordinates for the $N+1$ vertices and the one
curvature component for the cell giving a total of $2N+3$ unknowns. However, three of the
coordinates can be freely chosen (e.g., locate the origin at the central vertex and align
one axis with one of the legs) reducing the number of unknowns to $2N$. There are also $2N$
constraints provided by the known leg-lengths. Thus the coordinates and curvatures can be
found by solving the coupled system of equations given above (\ref{eqn:rncLsq}). One
objection to this method is that it can be computationally expensive to solve this system of
equations at each time step.

In the second approach, the Riemann curvatures are evolved using known evolution equations
for the Riemann curvatures. Morgan and Tian \cite{morgan:2007-01} show that\footnote{Here
the Riemann curvature is defined by $R^\alpha{}_{\mu\nu\beta}v^\beta =
v^\alpha{}_{\mu;\nu}-v^\alpha{}_{\nu;\mu}$ and the Ricci tensor by $R_{\mu\nu} =
R^\alpha{}_{\mu\nu\alpha}$}, under Ricci flow and in a local Riemann normal frame,
\begin{equation}
   \label{eqn:evolveRij}
   \frac{\partial R_{\alpha\beta}}{\partial t}
   =
   \nabla^2 R_{\alpha\beta}
   + 2 R^\mu{}_{\alpha\beta}{}^\nu R_{\mu\nu}
   - 2 R_{\alpha}{}^\mu R_{\beta\mu}
\end{equation}
(see equation 3.6 of \cite{morgan:2007-01}).

% see also Chow, Lu and Ni (2006) equation 2.3 and 2.37

For the case of a 2-dimensional metric, where $2R_{\mu\nu} = R g_{\mu\nu}$, it is easy to
show from (\ref{eqn:RicciFlowCoord}) and (\ref{eqn:evolveRij}) that
\begin{equation}
   \label{eqn:evolveR}
   \frac{\partial R}{\partial t}
   =
   \nabla^2 R + R^2
\end{equation}

\begin{figure}[ht]
\Figure{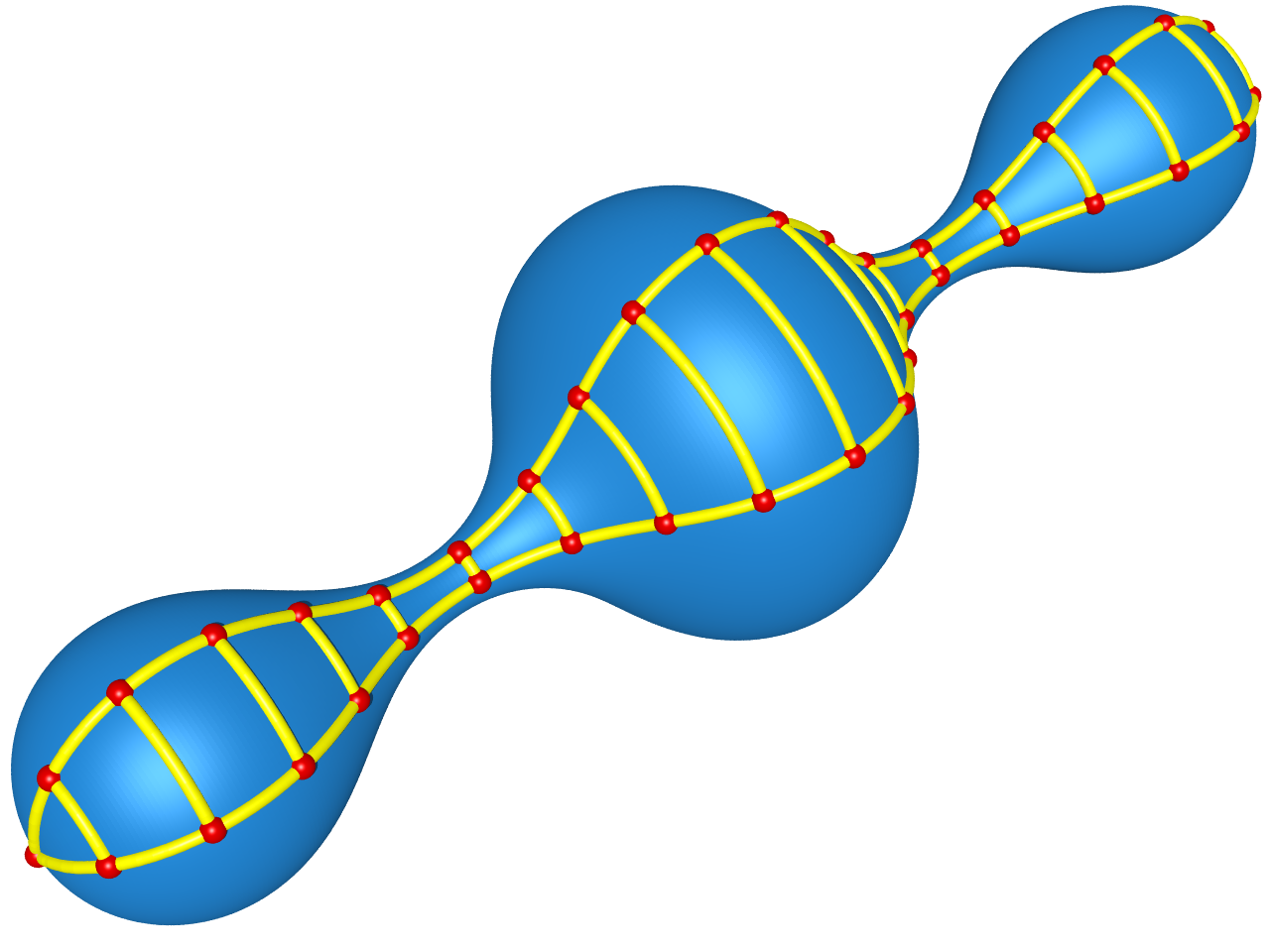}{0.85}
\vskip 10pt
\caption{\normalfont%
This is an example of the lattice used to represent the rotationally symmetric 2-geometries.
In the text this structure is often described as a \emph{ladder}. The ladder is bounded by
two geodesics that stretch from the north to south poles. The rungs of the ladder are short
legs that connect corresponding points on the rails of the ladder. A typical computational
cell is bounded by three consecutive rungs and the four segments of the rails of the ladder.
The region shared by two cells is bounded by a common pair of rungs and two segments of the
rails.}
\label{fig:ladder}
\end{figure}

A convenient lattice, for the 2-dimensional axisymmetric case, is shown in figure
(\ref{fig:ladder}). This has a ladder like structure in which each computational cell
consists of three consecutive rungs of the ladder. Pairs of cells overlap over a pairs of
rungs. A typical cell is shown in figure (\ref{fig:basiccell}) where the length of the rungs
are denoted by $L_x$ and the origin of the Riemann normal coordinates is chosen at the
centre of the rung and aligned with the $x$-axis pointing to the right and parallel to the
rung. It is possible to explicitly solve the leg-length equations (\ref{eqn:rncLsq}) (see
Appendix \ref{app:limitRNC} for full details) but for this particular lattice there is a
much quicker road to an equation linking the $L_x$ to the Riemann curvature. The axisymmetry
of the geometry shows that the rails of the ladder can be chosen to be geodesics that extend
from the north to south poles. Thus it is not surprising that a geodesic deviation equation
of the form
\begin{equation}
   \label{eqn:GeodesicDev}
   \frac{d^2L_{x}}{ds^2} = -\frac{1}{2} R L_{x}
\end{equation}
applies to this lattice\footnote{This equation can also be obtained as a limiting form of
the leg-length equations (\ref{eqn:rncLsq}) (as shown in Appendix \ref{app:limitRNC}).}.
This equation was used in the numerical simulations described below to compute $R$ given the
$L_x$ on the lattice.

The second method uses equation (\ref{eqn:evolveR}) to evolve the Ricci scalar. However this
requires values for $\nabla^2 R$ which in the local Riemann normal coordinates are given by
\begin{equation}
   \nabla^2 R = \frac{\partial^2 R}{\partial x^2} + \frac{\partial^2 R}{\partial y^2}
\end{equation}
At first glance it seems that a simple finite difference scheme could be used to estimate
the derivatives on the right hand side. This would be incorrect as it fails to take account
of the different coordinate frames being used in each cell. One way to deal with this issue
is to pick one cell and extend its coordinates into the neighbouring cells. This allows
coordinate transformations to be made so that data from neighbouring cells can be imported
into the chosen cell. At this point the finite difference approximation can be made. For our
choice of lattice this leads to (see Appendix \ref{app:nablaRNC})
\begin{equation}
   \label{eqn:flatLaplacian}
   \nabla^2 R
   = \frac{\partial^2 R}{\partial x^2} + \frac{\partial^2 R}{\partial y^2}
   = \frac{1}{L_x}\frac{dL_x}{ds}\frac{dR}{ds} + \frac{d^2R}{ds^2}
\end{equation}
where $s$ is the arc-length measured from the north and south poles and where $d/ds$ is the
ordinary derivative with respect to $s$.

In summary the evolution equations for this lattice are as follows. The leg-lengths are
evolved using (\ref{eqn:evolveLij}) while the curvatures are evolved using
\begin{equation}
   \label{eqn:rncEvolveR}
   \begin{aligned}
       R &= -\frac{2}{L_x}\frac{d^2L_{x}}{ds^2} &\qquad\text{method 1}&\\
       \frac{\partial R}{\partial t}
       &= R^2
       + \frac{1}{L_x}\frac{dL_x}{ds}\frac{dR}{ds}
       + \frac{d^2R}{ds^2} &\qquad\text{method 2}
   \end{aligned}
\end{equation}

% ============================================================================================
\section{Numerical formulation}
\label{sec:theNumerics}

\subsection{Initial data}

Rubinstein and Sinclair chose a two parameter family of initial data of the form
\begin{equation}
   \label{eqn:initialG}
   \begin{aligned}
      h(\rho) &= 1\\
      m(\rho) &= \left(\frac{\sin\rho + c_3\sin3\rho + c_5\sin5\rho}{1+3c_3+5c_5}\right)^2
   \end{aligned}
\end{equation}
for $0\le\rho\le\pi$ and where $c_3$ and $c_5$ are freely chosen constants. Note that this
choice of initial data is clearly consistent with the smoothness condition
(\ref{eqn:RSsmooth}).

The initial data for the Rubinstein and Sinclair model was obtained by sampling the above
metric functions on a uniform initial grid, namely, at $\rho_i = i\Delta\rho, i=0,1,2,\cdots
N-1$ with $\Delta\rho=\pi/N$. The same uniform grid was used to set the initial data for the
smooth lattice. Vertices were distributed uniformly from the north to south pole with labels
ranging from $0$ at the north pole to $N$ at the south pole. The length of the leg joining
vertex $i$ to $i+1$ (i.e., one segment of one rail of the ladder) is denoted by $L_{yi}$
while the length of the transverse leg passing through vertex $i$ (i.e., one of the rungs of
the ladder) is denoted by $L_{xi}$. Finally, the Ricci scalar at vertex $i$ is denoted by
$R_i$. The leg lengths $L_{xi}$ and $L_{yi}$ between the north pole and the equator were
computed using a numerical geodesic integrator while symmetry across the equatorial plane
was used to set the remaining leg lengths (Rubinstein and Sinclair do likewise for their
initial data).

\subsection{The poles}

Though the use of a numerical grid adapted to the symmetry reduces the computational
complexity it does introduce its own problems at the poles. This is clear from the form of
the evolution equations in which various terms are unbounded for grid points arbitrarily
close to the poles. Thus some care must be taken when evolving data at or close to the poles.

Rubinstein and Sinclair used equations (\ref{eqn:RSFlowB}) at the poles. They also noted
that local errors near the poles could violate the smoothness condition
(\ref{eqn:RSsmooth}). They dealt with this problem by replacing the evolved values of
$\sqrt{m}$ with $f(\rho)\sqrt{m(\rho)}$ with $f(\rho)$ a smooth function with values peaked
at the poles and rapidly decaying to 1 away from the poles. This function $f(\rho)$ was
specially crafted to preserve the smoothness condition (\ref{eqn:RSsmooth}).

\begin{figure}[ht]
\Figure{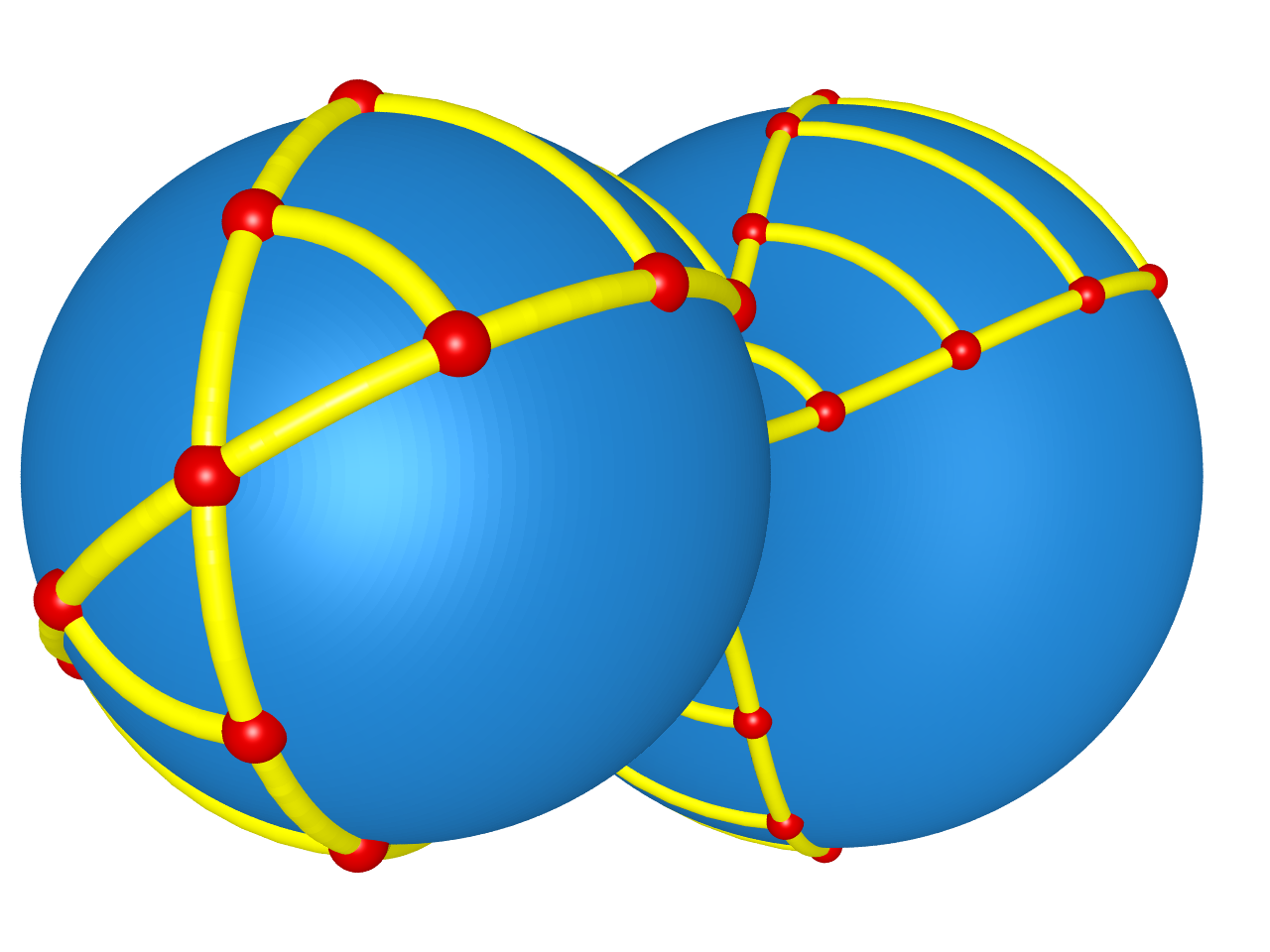}{0.85}
\vskip 10pt
\caption{\normalfont%
This is close-up view of the lattice showing clearly the extension of the lattice over the
poles. This extension is allowed by the rotational symmetry.}
\label{fig:over-poles}
\end{figure}

A different approach was employed for the smooth lattice equations. In this case the
symmetry of the underlying geometry was used to extend the data across the poles. For
example, the rails of the lattice are readily extended across the poles as shown in figure
(\ref{fig:over-poles}). Then at the north pole, $(L_{xi})_{i=-j} = - (L_{xi})_{i=j}$ while
$(L_{yi})_{i=-j} = (L_{yi})_{i=j+1}$ and $R_{-j} = R_{j}$ with a similar pattern applied at
the south pole. This makes it easy to compute symmetric finite difference approximations for
all non-singular terms in the neighbourhood of the poles. Let the vertices of the extended
lattice be labelled by $i$ with values $-m,-m+1,-m+2,\cdots 0,1,2,\cdots N,N+1,N+2,\cdots
N+m$ where $m$ is a small postive integer and where the poles have $i=0$ and $i=N$. For
method 1, the values of $R$ on the lattice were calculated as follows. Let $n<m$ be another
small postive integer. Then for $i=n,n+1,n+2\cdots N-n-2,N-n-1,N-n$ use equation
(\ref{eqn:rncEvolveR}) to compute $R$ while for the remaining values of $i$ use local
polynomial interpolation to fill in the remaining values for $R$. For method 2 there are two
choices available, either interpolate the time derivatives of $R$ befor a time step or
interpolate $R$ after a time step. The first choice failed badly for the double dumbbell but
work well for the 2-sphere and the single dumbbell. The second choice, to interpolate after
a time step, worked very well for all three initial data sets. Note that a fourth order
Runge-Kutta method (as used here) includes four time steps. The interpolation just described
was applied after each of the four time steps of the Runge-Kutta method.

Note that the rotational symmetry can be used to good effect in the interpolation. Consider
a smooth function $f(s)$ defined in a neighbourhood of the either pole and where $s$ is the
arc-length (along the geodesic) measured from the pole. Suppose $f(s)$ is an even function
of $s$ (such as for example $R$ and its time derivative). Then $f(s)$ can be expanded as a
power series
\begin{equation}
   f(s) = a_0 + a_2 s^2 + a_4 s^4 + \cdots
\end{equation}
for some set of coefficients $a_i,i=0,1,2,\cdots$ Let $f_i$ be the discrete value of $f(s)$
at $s=s_i$ for $i=n,n+1,n+2,\cdots m$. Then $f_j$ for $j=-n+1,-n+2,-n+3,\cdots 0,1,2,\cdots
n-1$ can be estimated by polynomial interpolation on the data $(s^2,f)_i,i=n,n+1,n+2\cdots
m$. This not only ensures smoothness across the poles but also creates higher order
estimates than would be obtained using the data $(s,f)_i$.

All of the results presented below were obtained using $n=2$ and $m=4$.

\subsection{Filtering}

In an attempt to minimise high frequency errors, Rubinstein and Sinclair filtered out high
frequency components in the gid values for $h(\rho)$ and $\sqrt{m(\rho)}$. They observed
that this reduced but did not cure the problem of numerical instabilities. No filtering of
this kind was used in the smooth lattice codes.

\subsection{Re-gridding}

Consider an initial lattice in which the vertices are uniformly distributed (i.e., $L_{yi}$
is constant over the range of $i$). At later times the vertices will no longer be uniformly
distributed with some vertices being drawn together while others will be pushed apart. This
is a natural outcome of the evolution under the Ricci flow. Unfortunately this can introduce
two problems. First, as some leg lengths shrink, the corresponding time step set by Courant
condition may become prohibitively small. Second, the ensuing irregular structure in the
lattice will lead to increasing truncation errors in the estimates of the first and second
derivatives required in equations (\ref{eqn:rncEvolveR}). These problems can be reduced by
periodic re-gridding of the lattice as follows. Create a new lattice, with uniformly
distributed $L_{yi}$, then use quadratic interpolation to produce new values for $L_{xi}$
and $R$ on the new lattice. This simple scheme proved to be the key step in obtaining long
term stable integrations. This observation was also noted by Rubinstein and Sinclair (using
a scheme almost identical to that used here though they refer to this as a reparametrization
of the data).

\subsection{Time step}

The Ricci flow equations are a form of heat equation and thus when using explicit forward
time integrators (such as a fourth-order Runge-Kutta for the smooth lattice or the FTCS used
by Rubinstein and Sinclair) a Courant like condition should be used to ensure stability of
the numerical solution. The Courant condition for a numerical heat equation is of the form
$\Delta t = C (\Delta x)^2$ where $C$ is the Courant factor, $\Delta t$ the time step and
$\Delta x$ a typical discretisation scale. For both smooth lattice methods it was found that
choosing $\Delta t = 0.1 \max_i L_{yi}$ worked very well with no signs of instabilties
throughout the simulation. It is interesting to note that Rubinstein and Sinclair chose a
fixed time step of $\Delta t = 0.0001$. This may explain the instabilities they observed
later in their evolutions.

\subsection{Embedding}

A natural question to ask is -- Can the 2-dimensional geometries generated by the Ricci flow
always be isometricly embedded in Euclidian $R^3$? In the general case, where no symmetries
apply, the answer is no. But for the case of a rotationally symmetric geometry Rubinstein
and Sinclair showed that an embedding is always possible. They went on to show that the
isometric surface in $R^3$ could be generated by rotating the curve described by
\begin{equation}
   \label{eqn:RSembedXY}
   \begin{aligned}
      x(\rho) &= \int_0^\rho \left(h(s)
               - \frac{1}{4m}\left(\frac{\partial m(s)}{\partial s}\right)^2\right)\>ds\\
      y(\rho) &= \sqrt{m(\rho)}
   \end{aligned}
\end{equation}
for $0\le\rho\le\pi$ around the $x$-axis (where $(x,y,z)$ are the usual Cartesian
coordinates covering $R^3$). Note that this curve passes through $(0,0)$. However for
aesthetic effect it is convenient to translate the curve along the $x$-axis so that $x(0) =
- x(\pi)$. For the initial data given above (\ref{eqn:initialG}) this produces a curve that
is reflection symmetric in the $y$-axis.

A similar construction should be possible starting form the lattice variables $L_x$ and
$L_y$. One approach would be to first extract the metric functions $h(\rho)$ and $m(\rho)$
from the lattice data $L_x$ and $L_y$ and to the then compute the generating curve using
equations (\ref{eqn:RSembedXY}). However, there is a simpler and more direct approach.
Consider two planes. The first plane $P$ is the $xy$-plane while the second plane $P'$
obtained by a small rotation of $P$ around the $x$-axis. The lattice, viewed as a ladder, is
then inserted between this pair of planes with each end of the ladder tied to the $x$-axis.
The Cartesian coordinates of each vertex can be computed from the $L_{xi}$ and $L_{yi}$ as
follows.

Let the $(x,y)_i$ be the coordinates of vertex $i$ in $P$. The coordinates, $(x,y)'_i$, of
the corresponding vertex in $P'$ will be obtained by a rotation of $(x,y)_i$ by an angle
$\alpha$, independent of $i$, around the $x$-axis. Thus
\begin{equation}
   \label{eqn:rotateXY}
   \begin{aligned}
      x'_i &= \phantom{-}x_i\cos\alpha + y_i\sin\alpha\\
      y'_i &= -x_i\sin\alpha + y_i\cos\alpha
   \end{aligned}
\end{equation}
Start by setting $(x,y)_0 = (0,0)$ and $(x,y)_1 = (0,L_{y0})$. Note that setting $x_0=x_1=0$
ensures that embedded surface is locally flat at the north pole. Then $(x,y)'_0 = (0,0)$ and
\begin{equation}
   \begin{aligned}
      x'_1 &= L_{y0}\sin\alpha\\
      y'_1 &= L_{y0}\cos\alpha
   \end{aligned}
\end{equation}
But the leg joining $(x,y)_1$ to $(x,y)'_1$ (i.e., a rung of the ladder) has length $L_{x1}$
and thus using the Euclidean metric of $R^3$ leads to
\begin{equation}
   L^2_{x1} = (x'_1-x_1)^2 + (y'_1-y_1)^2 = 2L^2_{z0}(1-\cos\alpha)
\end{equation}
which allows $\alpha$ to be computed from $L_{x1}$ and $L_{y0}$. The remaining coordinates
are computed in a similar manner. Suppose the coordinates for vertices $0,1,2,\cdots i-1$
have been computed. Then the coordinates for the next vertex $i$ are obtained by solving the
coupled pair of equations
\begin{equation}
   \label{eqn:rncEmbedXY}
   \begin{aligned}
      L^2_{zi-1} &= (x_i-x_{i-1})^2 + (y_i-y_{i-1})^2\\
      L^2_{xi}   &= (x'_i-x_i)^2 + (y'_i-y_i)^2
   \end{aligned}
\end{equation}
for $(x,y)_i$. Once all of the coordinates have been computed the curve is translated along
the $x$-axis, to centre the curve, by the replacement $x_i \mapsto x_i - x_e$ where $x_e$ is
the (original) $x$-coordinate of the vertex on the equator.

% ============================================================================================
\section{Results}
\label{sec:results}

Results for three distinct initial datasets, two that match those of Rubinstein and Sinclair
and a third for a unit 2-sphere, will be presented.

The first dataset, which has the look of a single dumbbell, uses $c_3=0.766,c_5=-0.091$
while the second dataset, a double dumbbell, has $c_3=0.021,c_5=0.598$ and the third is a
unit 2-sphere with $c_3=c_5=0$.

The slrf codes were run until the time step had been reduced by a pre-determined factor
(200 for the 2-sphere and 400 for the dumbells). The ricci-rot code ran until it detected
significant numerical errors (beyond this point the code would crash). The run time for the
ricci-rot codes was always shorter than that for the slrf-codes.

\begin{figure}[ht]
\Figure{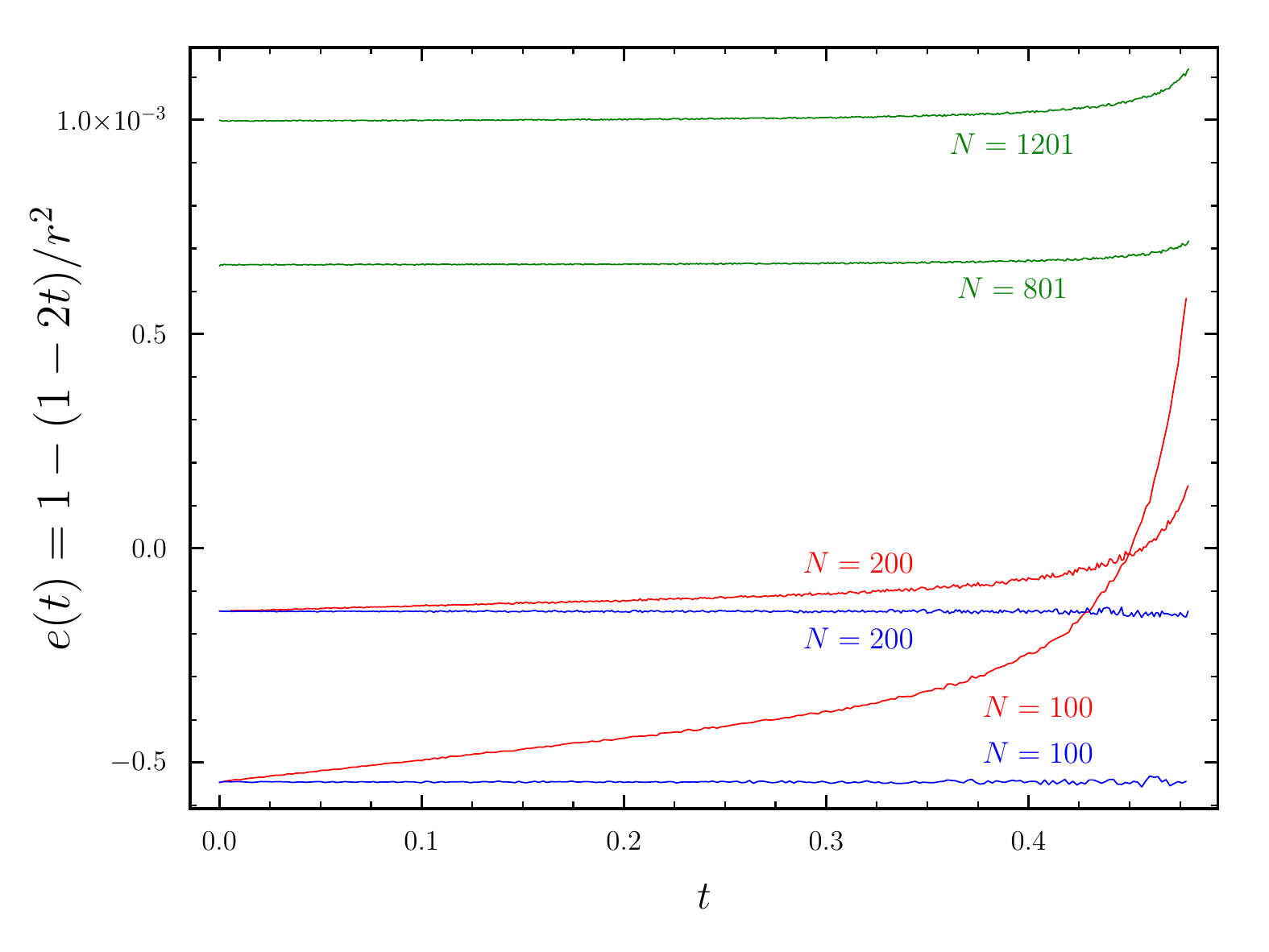}{0.95}
\caption{\normalfont%
These plots display the history of the fractional error in $r^2(t)$ for the Ricci flow of a
unit 2-sphere. The green curves are those for the ricci-rot code, while the remaining curves
are for the smooth lattice methods (red for the first method, blue for the second). The
noise in the curves is due to the regridding (for the smooth lattice codes) and the
reparametrization (for the ricci-rot code). The upward rise in the errors for late times is
most likely due to to increased truncation errors.}
\label{fig:2sphere12}
\end{figure}

For initial data based on a unit 2-sphere it is a easy to show that the subsequent evolution
continues to be a 2-sphere with a radius $r(t)$ that evolves according to
\begin{equation}
   r^2(t) = 1 - 2 t
\end{equation}
This not only shows that $r\rightarrow0$ as $t\rightarrow 1/2$ but it also provides a very
simple test of each of the three computer codes, the Rubinstein and Sinclair code (which
they named ricci-rot) and the two smooth lattice codes (which will be referred to as slrf-v1
and slrf-v2 corresponding to the two smooth lattice methods). Figure (\ref{fig:2sphere12})
displays the the fractional error in $r^2(t)$, defined by
\begin{equation}
   e(t) = \frac{r^2(t) -1 + 2t}{r^2(t)}
\end{equation}
as a function of time. It shows clearly that, for this initial data, the three codes produce
very good results. The upward trend in the curves near $t=0.5$ is most likely a numerical
artefact due to the rapidly increasing curvature ($R\approx 500$ at $t=0.49$). The small
noise in the curves are due to the re-gridding (slrf-v1,v2) and filtering (ricci-rot).

\begin{figure}[ht]
\Figure{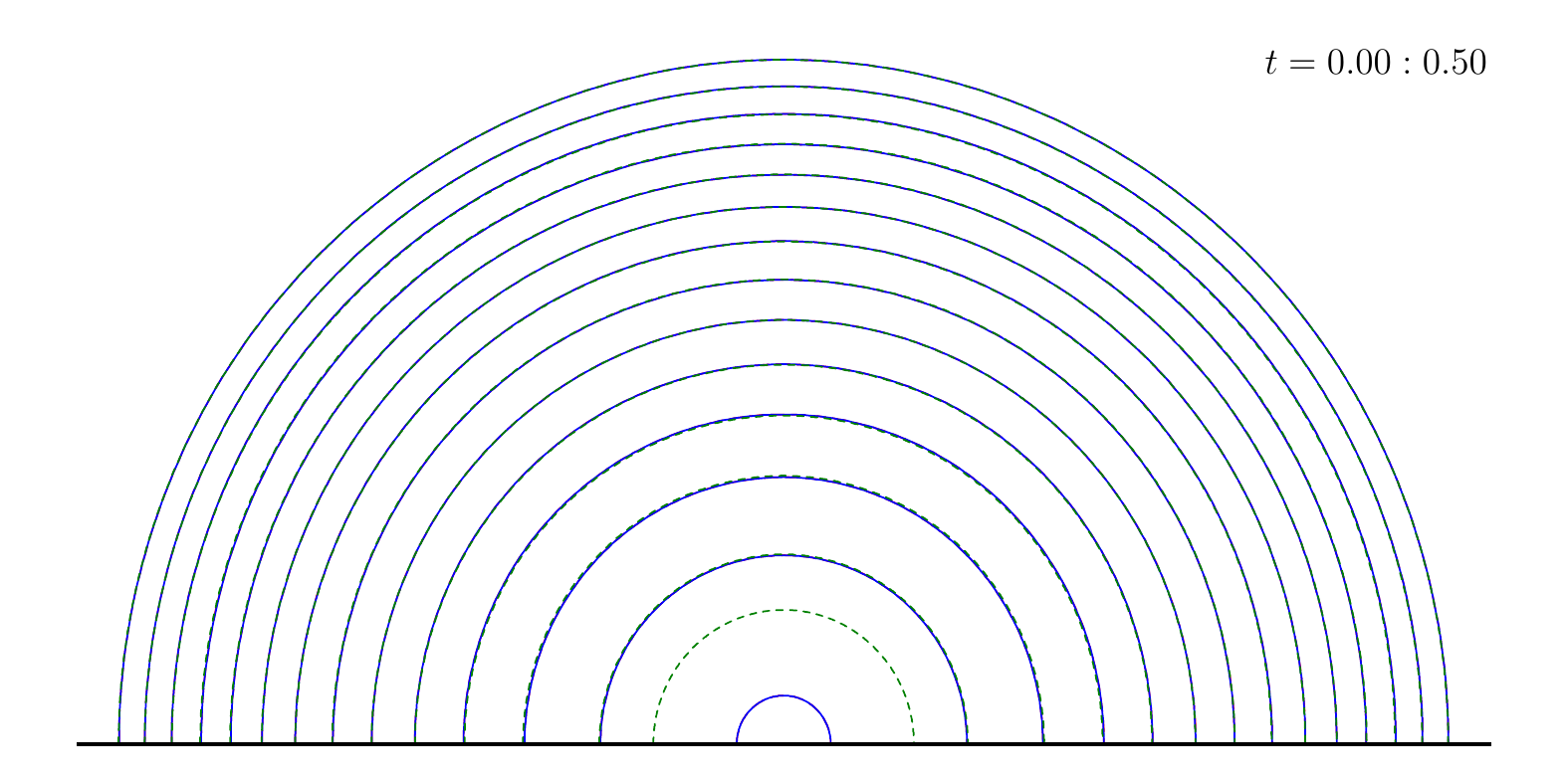}{0.90}
\Figure{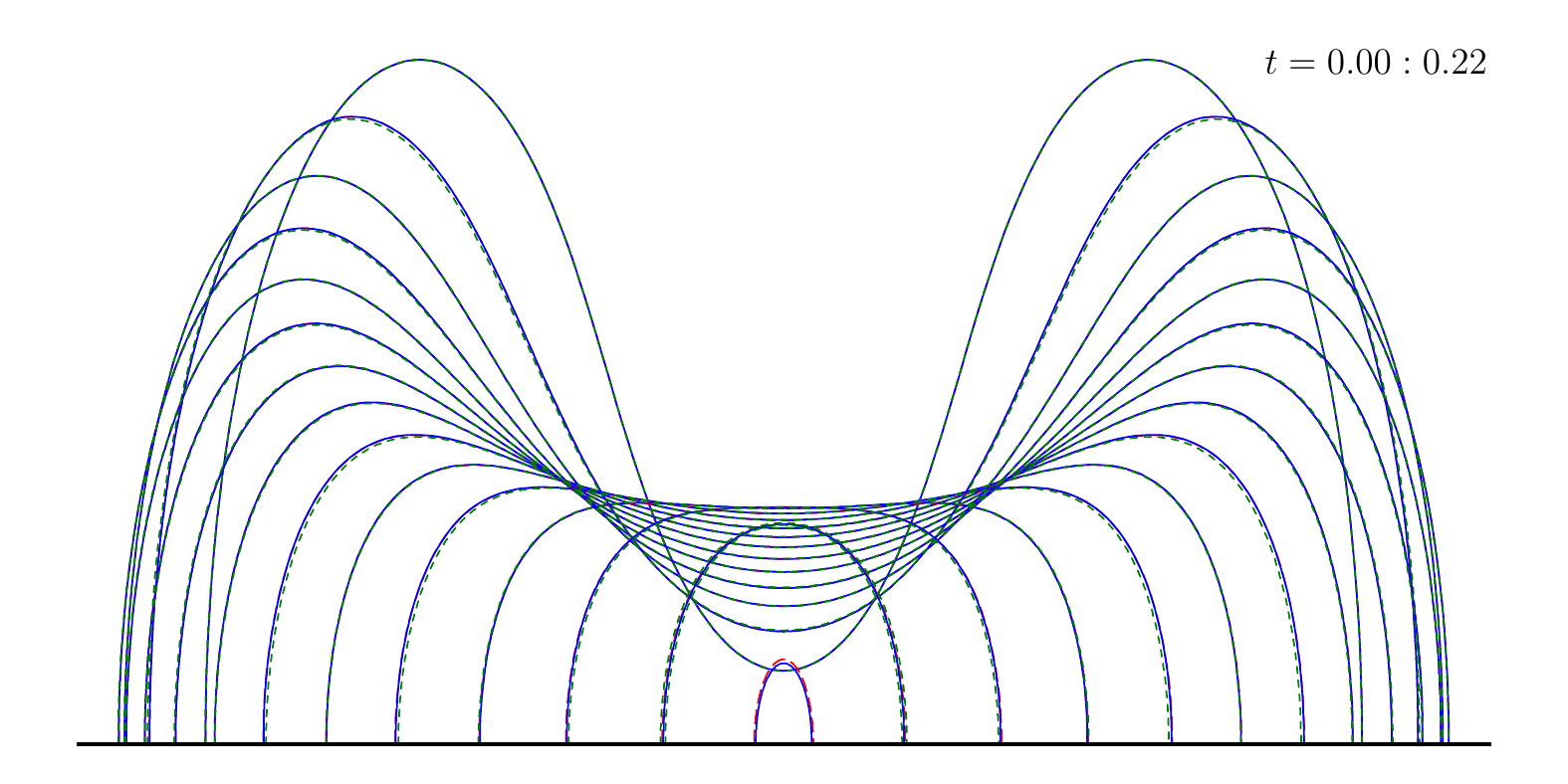}{0.90}
\caption{\normalfont%
Each curve in each of these plots is a generating curve for the isometric embedding of the
2-geometry in $E^3$. The full 2-geometry would be obtained by rotating the curve around the
horizontal axis. The upper diagram is for the 2-sphere and shows clearly that the geometry
remains spherical throughout the evolution. The lower diagram is for the single dumbbell
initial data. This shows a transition from an initial geometry with both postive and
negative $R$ through to later stages where $R$ is strictly positive and increasing (as
expected). Note that each curve here is actually three curves, one for each of the three
methods. That they appear as a single curve is strong evidence of the correctness of the
numerical codes.}
\label{fig:2sphere1xdbell}
\end{figure}

\begin{figure}[ht]
\Figure{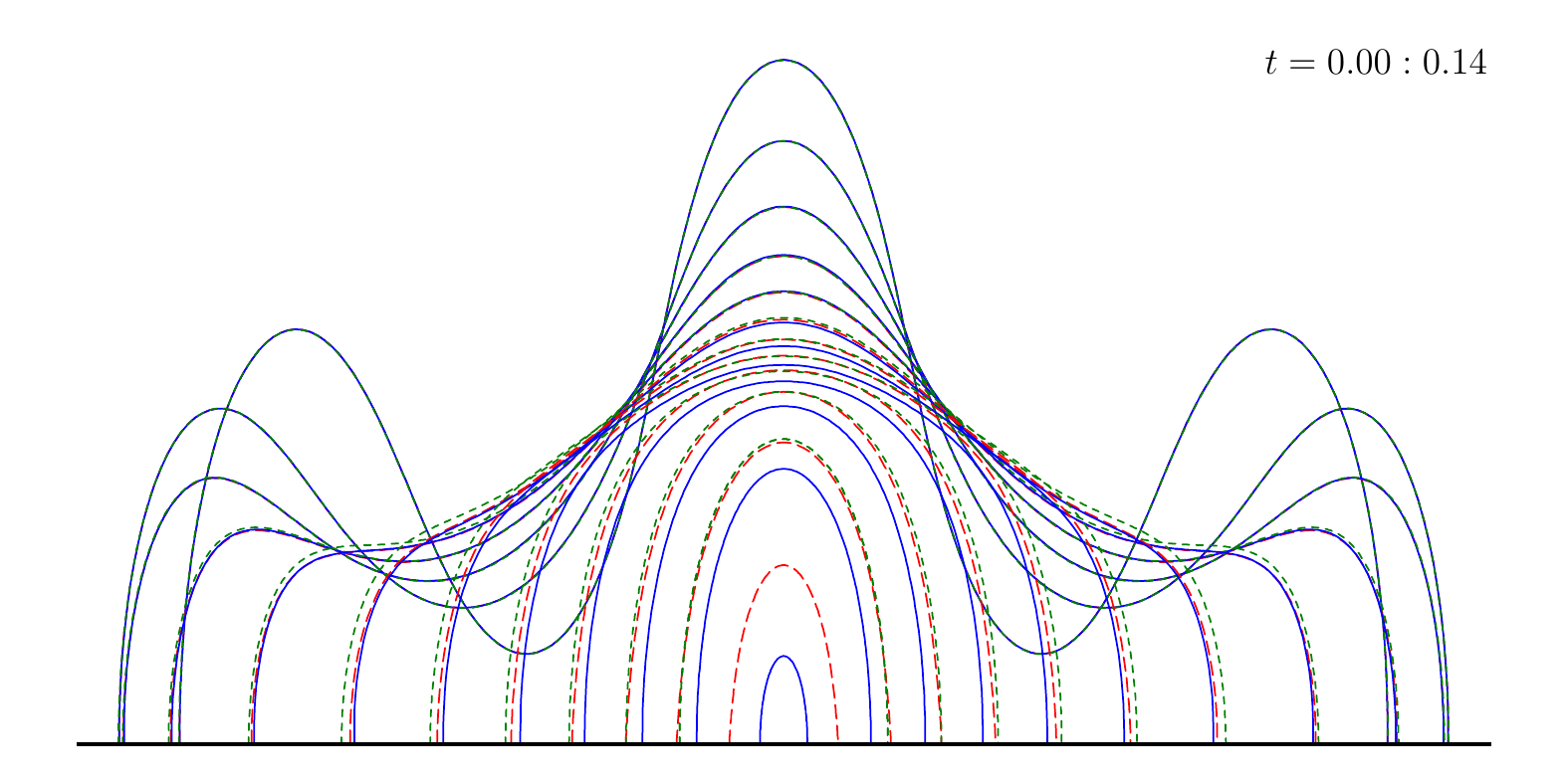}{0.90}
\Figure{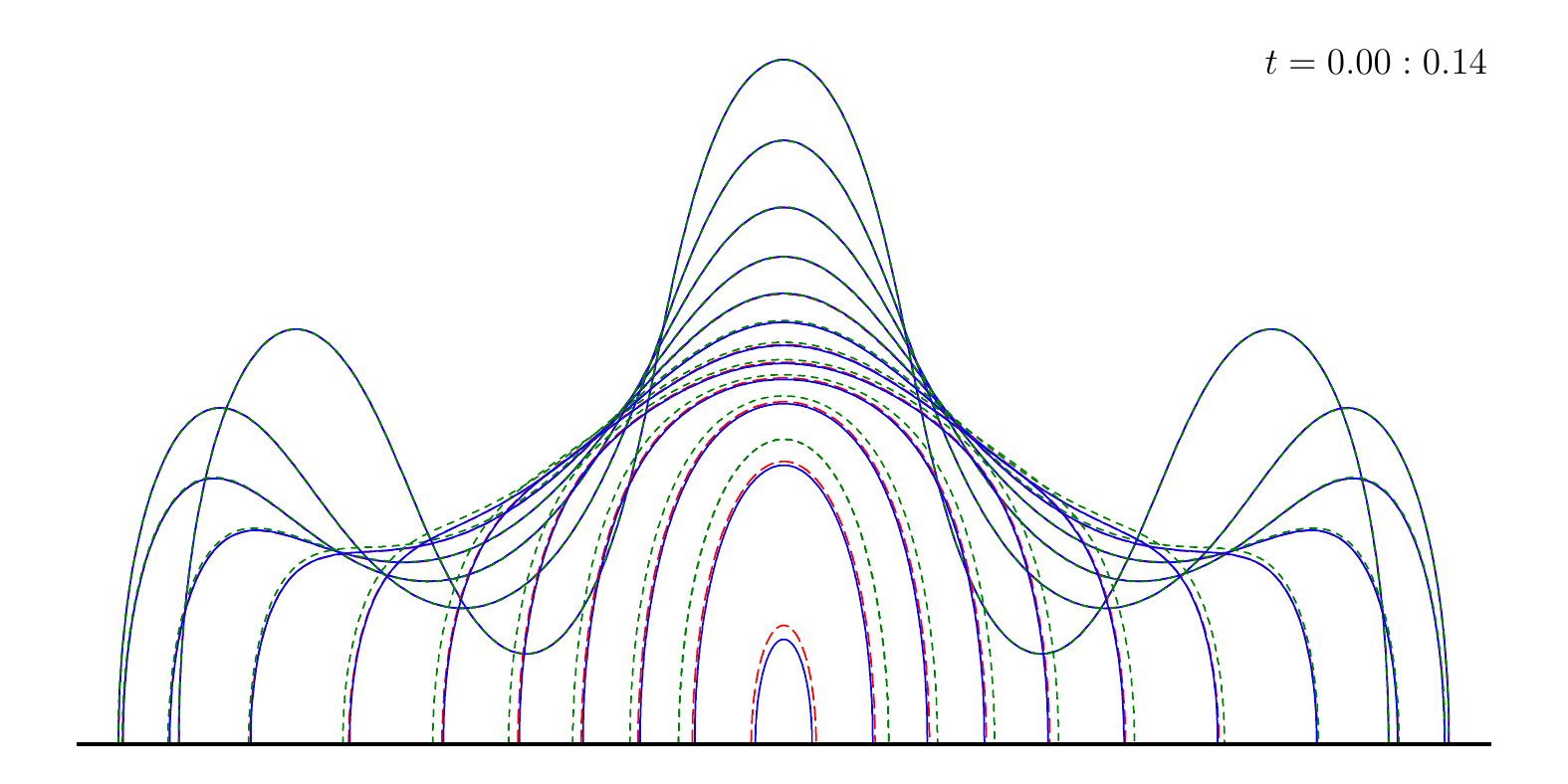}{0.90}
\caption{\normalfont%
This shows the generating curves for the double dumbbell initial data. Unlike the previous
figure, here it is possible to discern small differences between the three methods for late
times. The upper figure was run with a coarse resolution ($N=100$ for the slrf codes,
$N=801$ for the ricci-rot code) while the lower figure uses $N=200$ and $N=1201$
respectively. This increased resolution has brought the three curves closer together and
thus the differences are most likely due to limited resolution rather than any error in the
code.}
\label{fig:2xdbell}
\end{figure}

Figures (\ref{fig:2sphere1xdbell},\ref{fig:2xdbell}) display a history of the embeddings for
each of the three models with each figure showing results from all three codes (short dashes
for ricci-rot, a longer dashes for slrf-v1 and a solid line for slrf-v2). The ricci-rot code
generally terminated much earlier than the slrf codes\footnote{This could probably be
improved by allowing a variable time step in the ricci-rot code.} and this can be seen in the
upper panel of figure (\ref{fig:2sphere1xdbell}) which shows a lone dashed curve this being
the final curve computed by the ricci-rot code. Similar lone curves can be seen in figure
(\ref{fig:2xdbell}).

For the 2-sphere and single dumbbell all three codes gave almost identical results (the
three curves appearing almost as a single curve) but for the double dumbbell there is a small
difference between the codes at later times. This difference is almost certainly due to
discretisation errors. This claim was tested by running the codes with an increased number
of grid points (from 100 to 200 for slrf-v1,v2 and from 801 to 1201 for
ricci-rot\footnote{The ricci-rot code was unable to evolve the double dumbbell initial data
to $t=0.11$ for 1601 grid points.}). The results are shown in figure (\ref{fig:2xdbell})
which shows improved agreement between all three codes.

The fact that all three codes produce almost identical results should not be understated.
The codes employ fundamentally different algorithms, they were written by different people
and in different computer languages (ricci-rot in C and slrf-v1,v2 in Ada). That these codes
agree as well as they do is a very strong indication that they are giving correct results.

\section{Discussion}
\label{sec:dicuss}
Though the results presented here for the smooth lattice method are encouraging much work
still remains to be done. For example, can the method be used for non-symmetric geometries
in two or more dimensions? Another line of investigation would to be to compare the relative
merits of the two smooth lattice methods presented here. It might be argued that in the
absence of global geodesics the first method would be ruled out as the geodesic deviation
equation could not be readily adapted to the lattice. However, it would still be possible to
extract estimates for $R$, without resort to the geodesic equation, as shown in appendix
(\ref{app:limitRNC}). These and other questions will be explored in a later paper.

The clear separation between the metric and topology on the lattice seems well suited to the
study of Ricci flow in 3 dimensions where complex behaviours are known to develop.

\appendix

% ============================================================================================
\section{The geodesic deviation equation}
\label{app:limitRNC}

In section (\ref{sec:SLRF}) it was claimed that the geodesic deviation equation
(\ref{eqn:GeodesicDev}) can be obtained from the leg-length equation (\ref{eqn:rncLsq})
provided that the rails of the ladder (see figure (\ref{fig:ladder})) are be chosen to be
geodesics of the 2-geometry. The purpose of this appendix is to fill in the details of that
claim.

\begin{figure}[ht]
\Figure{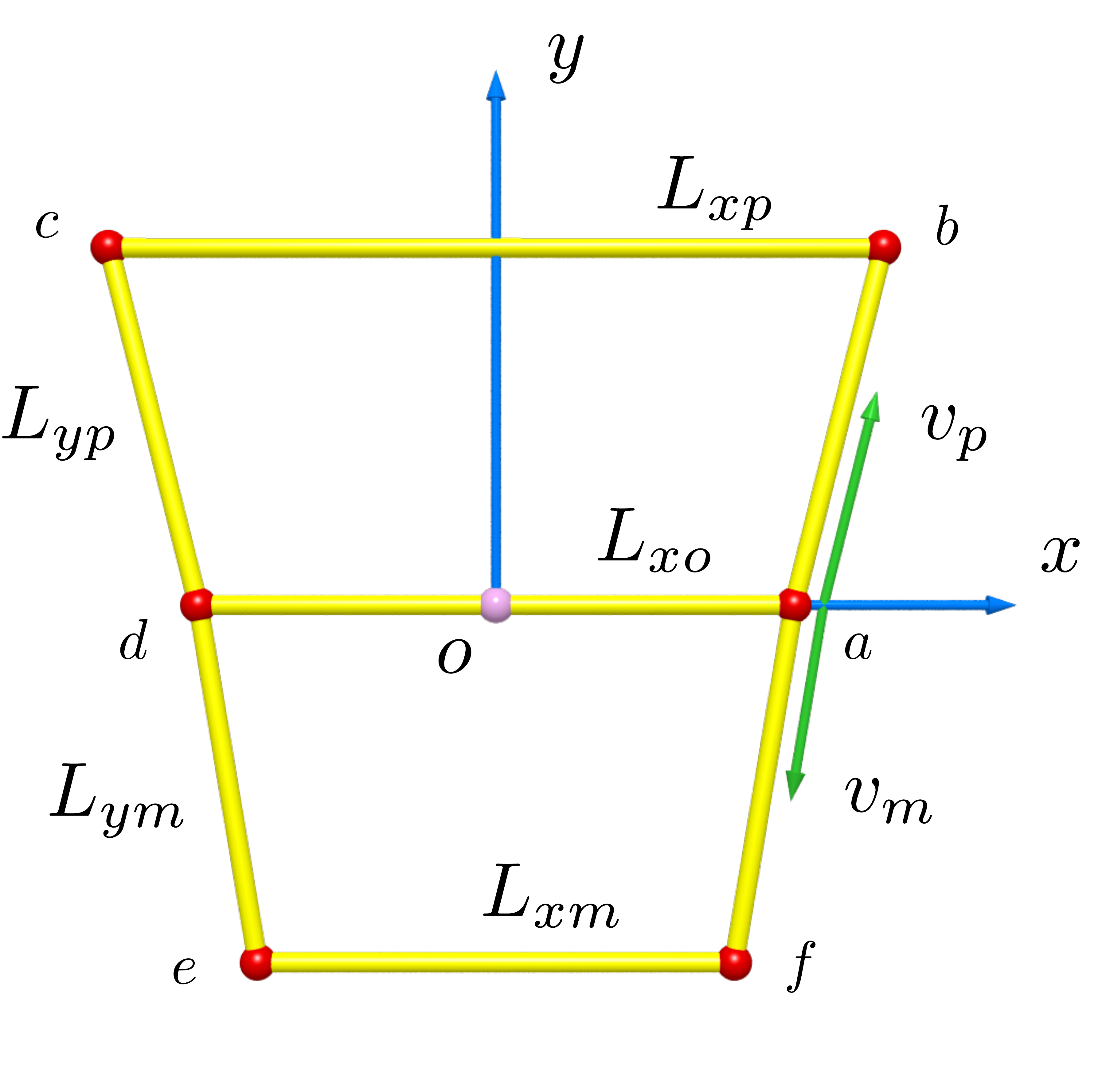}{0.90}
\caption{\normalfont%
A typical computational cell in the ladder lattice. The horizontal legs are the rungs of the
ladder while the remaining legs are segments of the rails of the ladder. In the text the
typical leg-lengths are denoted by $L_x$ and $L_y$ however in this diagram extra sub-scripts
$o,p$ and $m$ are used to distinguish the various legs. The origin of the Riemann normal
coordinate frame is located at the vertex $o$. The vectors $v_p$ and $v_m$ are the unit
tangent vectors to the geodesic segments. The vertices are labeled clockwise from $a$ to
$f$. The coordinates for these vertices in terms of the leg-lengths are derived in Appendix
\ref{app:limitRNC}.}
\label{fig:basiccell}
\end{figure}

A single computational cell is shown in figure (\ref{fig:basiccell}). This consists of seven
vertices and eight legs. The Riemann normal coordinates are chosen so that the origin is
located at the centre of leg $(ad)$ and the $x$-axis is aligned to the same leg. The assumed
symmetries in the 2-geometry ensures that the coordinates for each of the seven vertices can
be chosen as shown in table (\ref{tbl:RNCCoords}).

\bgroup
\def\H{\vrule height 14pt depth  7pt width 0pt}
\def\m{\vrule height  0pt depth  8pt width 0pt}
\def\M{\vrule height 15pt depth 10pt width 0pt}
\def\A#1{\hbox to 40pt{\hfill$#1$}}
\def\B#1{\hbox to 20pt{\hfill$#1$}}
\def\C#1{\hbox to 12pt{\hfill$#1$}}
\def\D#1{\hbox to 10pt{\hfill$#1$}}
\begin{table}[t]
\begin{center}
\begin{tabular}{ccccccc}
\hline
% \H&Vertex&&\multicolumn{3}{c}{Coordinates}&\\
% \hline
\M&$a=(\A{ L_{xo}/2},\D{0})$&&$b=(\B{ x_b},\C{y_b})$&&$f=(\B{ x_f},\C{y_f})$&    \\
\m&$d=(\A{-L_{xo}/2},\D{0})$&&$c=(\B{-x_b},\C{y_b})$&&$e=(\B{-x_f},\C{y_f})$    \\
\hline
\end{tabular}
\end{center}
\caption{The Riemann normal coordinates of the 6 of the 7 vertices in figure
(\ref{fig:basiccell}). The remaining vertex $o$ has coordinates $(0,0)$.}
\label{tbl:RNCCoords}
\end{table}
\egroup

With this choice of coordinates it is a simple matter to write out in full the leg-length
equations (\ref{eqn:rncLsq}). This leads to
\begin{equation}
   \label{eqn:theLsqs}
   \begin{aligned}
      L^2_{px} &= L^2_{bc} = 4x^2_b - \frac{2}{3}(x_b y_b)^2 R\\
      L^2_{mx} &= L^2_{ef} = 4x^2_f - \frac{2}{3}(x_f y_f)^2 R\\
      (\Delta s_p)^2 &= L^2_{ab} = \frac{1}{4}(2x_b-L_{xo})^2
                                      + y^2_b - \frac{1}{6}(L_{xo}y_b)^2 R\\[5pt]
      (\Delta s_m)^2 &= L^2_{af} = \frac{1}{4}(2x_f-L_{xo})^2
                                      + y^2_f - \frac{1}{6}(L_{xo}y_f)^2 R
   \end{aligned}
\end{equation}
where $\Delta s_p$ and $\Delta s_m$ are defined by $\Delta s_p = L_{yp}$, $\Delta s_m =
L_{ym}$ and $R=2R_{xyxy}$ is the Ricci scalar. This is a non-linear system of four equations
for five unknowns, namely, the curvature $R$ and the four coordinates $(x_b,y_b)$ and
$(x_f,y_f)$. Put aside, for the moment, the obvious problem that an extra equation is
required to fully determine all five unknowns. Then the above equations can be viewed as a
set of four equations for the four unknown coordinates. Consider now the progression towards
the continuum limit where the computational cells will be much smaller than the length scale
associated with the curvature (e.g., $L^2_{ab}\ll1/R$). It follows that the curvature terms
in the above equations will be small compared to the leading terms and thus it should be
possible to express the coordinates as a power series in $R$. Thus put
\begin{equation}
   \begin{aligned}
      x_b & = x_{0b} + x_{1b} R + \BigO{R^2}\\
      x_f & = x_{0f} + x_{1f} R + \BigO{R^2}
   \end{aligned}
\end{equation}
Similar expansions could be made for $y_b$ and $y_f$. However, equations (\ref{eqn:theLsqs})
contain $y^2_b$ and $y^2_f$ but not $y^1_b$ nor $y^1_f$. Thus it is simpler to use a power
series for $y^2$ rather than for $y$, that is
\begin{equation}
   \begin{aligned}
      y^2_b & = w_{0b} + w_{1b} R + \BigO{R^2}\\
      y^2_f & = w_{0f} + w_{1f} R + \BigO{R^2}
   \end{aligned}
\end{equation}
Substituting these into equations (\ref{eqn:theLsqs}) and then following the standard
procedure of expansion and equating terms leads to
\begin{equation}
   \begin{aligned}
      0 &= -L^2_{px} + 4 x^2_{0b}\\
      0 &= -L^2_{mx} + 4 x^2_{0f}\\
      0 &= -4(\Delta s_p)^2 + 4w_{0b} + (2x_{0b}-L_{xo})^2\\
      0 &= -4(\Delta s_m)^2 + 4w_{0f} + (2x_{0f}-L_{xo})^2
   \end{aligned}
\end{equation}
for the $R^0$ terms and
\begin{equation}
   \begin{aligned}
      0 &= 12x_{1b}x_{0b} - x^2_{0b}w_{0b}\\
      0 &= 12x_{1f}x_{0f} - x^2_{0f}w_{0f}\\
      0 &= 6w_{1b} + 6(2x_{0b}-L_{xo})x_{1b} - w_{0b}L^2_{ox}\\
      0 &= 6w_{1f} + 6(2x_{0f}-L_{xo})x_{1f} - w_{0f}L^2_{ox}
   \end{aligned}
\end{equation}
for the $R^1$ terms. This set of equations are easily solved leading to
\begin{equation}
   \begin{aligned}
      96 x_b &= 48L_{xp} - \left((L_{xo}-L_{xp})^2-4(\Delta s_p)^2\right)RL_{xp}\\
      96 x_f &= 48L_{xm} - \left((L_{xo}-L_{xm})^2-4(\Delta s_m)^2\right)RL_{xm}\\
      96 y^2_b &=48\left(4(\Delta s_p)^2-(L_{xo}-L_{xp})^2\right)\\
               &- \left((L_{xo}-L_{xp})^2-4(
                          \Delta s_p)^2\right)(4L^2_{ox}+L_{xp}L_{xo}-L^2_{px})R\\[3pt]
      96 y^2_f &= 48\left(4(\Delta s_m)^2-(L_{xo}-L_{xm})^2\right)\\
               &- \left((L_{xo}-L_{xm})^2-4(
                          \Delta s_m)^2\right)(4L^2_{ox}+L_{xm}L_{xo}-L^2_{mx})R
   \end{aligned}
\end{equation}

There remains of course the issue of the missing equation. Since there are no more
leg-length equations to invoke that final equation must come from some constraint to remove
any remaining freedoms in the structure of the lattice. There are considerable freedoms in
choosing where to locate the vertices on the 2-geometry. An obvious choice is to require
that the edges that comprise the rails of the ladder form a global geodesic (i.e., a single
geodesic that stretches form the north to south pole). This can be achieved by requiring the
tangent vectors along the rails to be continuous from one segment to the next (e.g., from
leg $(fa)$ to leg $(ab)$). Continuity of the tangent vector at vertex $a$ in figure
(\ref{fig:basiccell}) thus requires
\begin{equation}
   0 = v^\mu_p + v^\mu_m
\end{equation}
There are actually two equations here, one for each component of the vectors. However, the
$y$ component yields a trivial equation in the continuum limit which leaves just the $x$
component as the required final equation. Using equation (\ref{eqn:rncGeodesicTangent}) to
construct the unit tangent vectors $v_p$ and $v_m$ and then setting $v^x_p+v^x_m$ to zero
leads to the following equation
\begin{equation}
   \begin{aligned}
   0 &= \frac{1}{2}\frac{L_{xp}-L_{xo}}{\Delta s_p}
      + \frac{1}{2}\frac{L_{xm}-L_{xo}}{\Delta s_m}\\
     &- \frac{(2L_{xo}+L_{xp})\left((L_{xo}-L_{xp})^2-4(\Delta s_p)^2\right)R}{96\Delta  s_p}\\
     &- \frac{(2L_{xo}+L_{xm})\left((L_{xo}-L_{xm})^2-4(\Delta s_m)^2\right)R}{96\Delta s_m}
   \end{aligned}
\end{equation}
It is not hard to see that this is a non-uniform finite difference approximation for the
following differential equation
\begin{equation}
   0 = \frac{d^2L}{ds^2} + \frac{1}{2} R L - \frac{1}{8} R L \left(\frac{dL}{ds}\right)^2
\end{equation}
where now $L=L(s)$. This equation can be reduced to the geodesic deviation equation in the
case where $L$ is taken to be very small relative to the curvature scales.

% ============================================================================================
\section{The Laplacian}
\label{app:nablaRNC}

In this appendix the details of the calculations leading to equation
(\ref{eqn:flatLaplacian}) will be presented.

As noted in the text (section \ref{sec:SLRF}), a correct calculation of $\nabla^2R$ will
require proper attention to the coordinate transformations required when sharing data
between neighbouring frames. The calculation of $\nabla^2 R$ will proceed by way of a finite
difference method but not before the relevant data has been imported from the neighbouring
frames. The first job then is to build these transformations.

\begin{figure}[ht]
\Figure{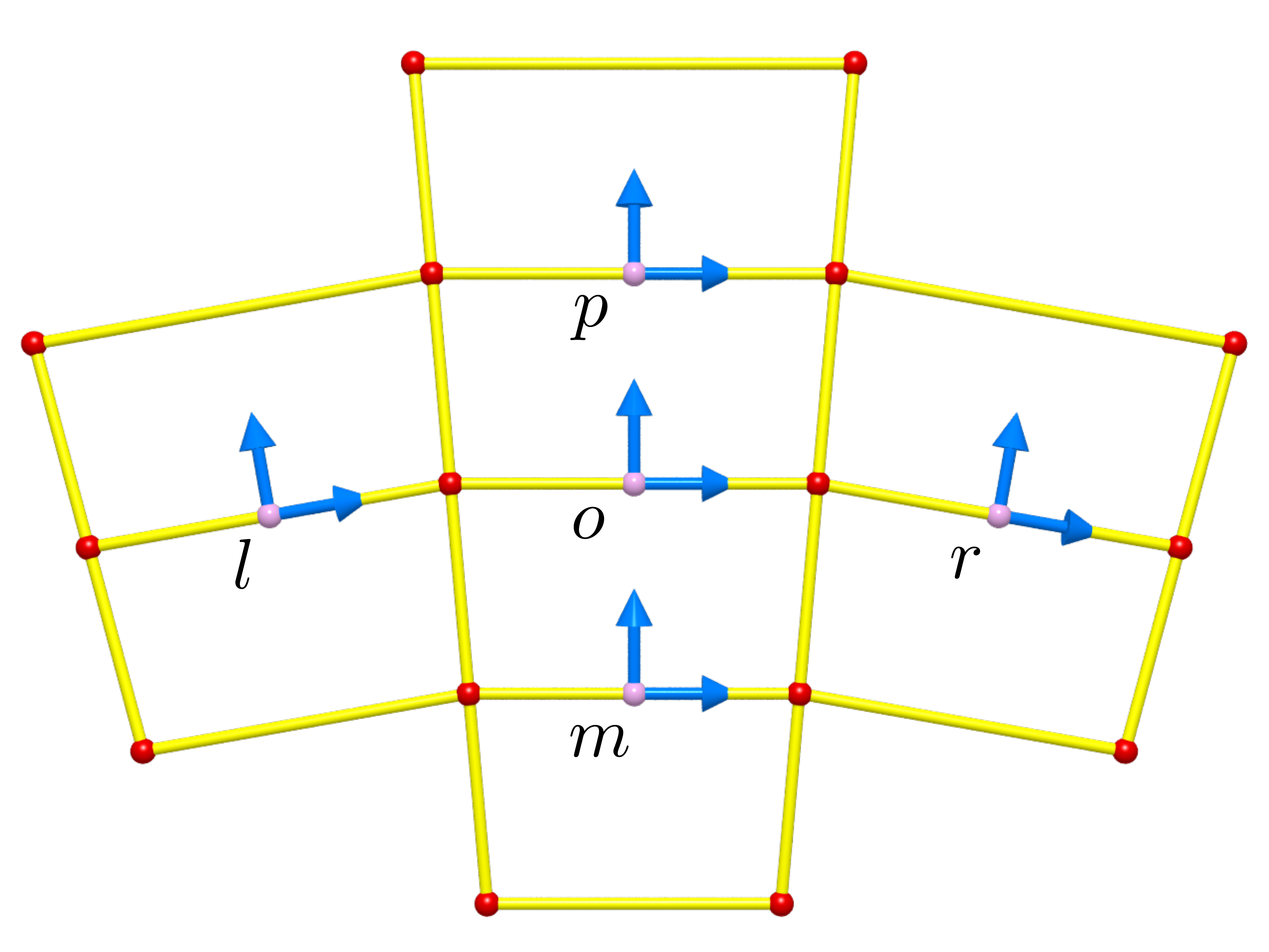}{0.90}
\caption{\normalfont%
This diagram shows a single computational cell labelled $o$ and its for four neighbouring
cells. The data in the cells labelled $l$ (for \emph{left}) and $r$ (for \emph{right}) are
identical to that in cell $o$ (due to rotationally symmetry). The arrows in each cell
represent the axes of the Riemann normal coordinate frames. The rotation angle from cell $o$
to cell $r$ can be computed using standard flat space Euclidian geometry (the rotation
matrix is applied in Appendix \ref{app:nablaRNC} only to the Riemann curvatures and thus to
leading order in the curvatures, the Euclidian approximation is appropriate). Note that the
overlap between cells $r$ and $o$ contains just two legs (likewise for cells $l$ and $o$).
This is not sufficient to fully determine the coordinate transformation between the
respective frames. This problem is resolved by using the known rotational symmetry.}
\label{fig:5cells}
\end{figure}

Figure (\ref{fig:5cells}) shows a set of 5 computational cells. Cells $m,o$ and $p$ are
three consecutive cells on the lattice while cells $l$ and $r$ are ghost cells created as
clones of cell $o$. The data in the ghost cells $l$ and $r$ are, by axisymmetry, identical
to the data in cell $o$. However the basis vectors of their coordinate frames are not
aligned to each other. It is easy to see that there is a simple rotation that maps the basis
vectors from one cell to the other. Consider now a vector field $v$ on the 2-dimensional
manifold. Let ${\bar p}$ denote the coordinate frame associated with the cell $p$. Let
$v^\mu_{a{\bar b}}$ denote the components of $v$ at the point $a$ in the frame ${\bar b}$.
Suppose that $v$ is a vector that respects the axisymmetry of the 2-geometry, then
\begin{equation}
   v^\mu_{l{\bar l}} = v^\mu_{o{\bar o}} = v^\mu_{r{\bar r}}
\end{equation}
Now at the point $l$ the basis vectors for frame ${\bar l}$ are a rotation of those for
${\bar o}$. Thus
\begin{equation}
   \begin{aligned}
      v^x_{l{\bar o}} &= v^x_{l{\bar l}}\cos\theta - v^y_{l{\bar l}}\sin\theta\\
      v^y_{l{\bar o}} &= v^y_{l{\bar l}}\cos\theta + v^x_{l{\bar l}}\sin\theta
   \end{aligned}
\end{equation}
for some angle $\theta$. This is one step in the process of importing data into frame ${\bar
o}$. The same argument can be applied for right hand frame ${\bar r}$. This leads to
\begin{equation}
   \begin{aligned}
      v^x_{r{\bar o}} &= v^x_{r{\bar r}}\cos\theta + v^y_{r{\bar r}}\sin\theta\\
      v^y_{r{\bar o}} &= v^y_{r{\bar r}}\cos\theta - v^x_{r{\bar r}}\sin\theta
   \end{aligned}
\end{equation}
Clearly no rotations are required for the frames ${\bar m}$, ${\bar o}$ and ${\bar f}$ and
thus
\begin{equation}
   \begin{aligned}
      v^x_{p{\bar o}} &= v^x_{p{\bar p}}\\
      v^y_{p{\bar o}} &= v^y_{p{\bar p}}\\
      v^x_{m{\bar o}} &= v^x_{m{\bar m}}\\
      v^y_{m{\bar o}} &= v^y_{m{\bar m}}
   \end{aligned}
\end{equation}
The components of the vector $v$ are now known at all five points in the one frame, ${\bar
o}$. Thus the derivatives at $o$ can be estimated using a finite difference method. For
example
\begin{equation}
   \frac{\partial v^x}{\partial x}
       = \frac{v^x_{r{\bar o}}-v^x_{l{\bar o}}}{L_{ox}}
       = \frac{v^y_{r{\bar r}}\sin\theta - v^y_{l{\bar l}}\sin\theta}{L_{ox}}
\end{equation}
But $v^y_{r{\bar r}} = v^y_{l{\bar l}} = v^y_{o{\bar o}}$ and thus
\begin{equation}
   \frac{\partial v^x}{\partial x} = v^y_{o{\bar o}}\frac{2\sin\theta}{L_{ox}}
\end{equation}
and to leading order in the lattice scale $2\sin\theta = (L_{px}-L_{ox})/\Delta s = dL_x/ds$
thus
\begin{equation}
   \frac{\partial v^x}{\partial x} = v^y \frac{1}{L_x}\frac{dL_x}{ds}
\end{equation}
where $v^x_{o{\bar o}}$ has been abbreviated to $v^x$ and $L_{ox}$ to $L_x$. Since the
transformations between the frames ${\bar m}$, ${\bar o}$ and ${\bar f}$ are trivial there
is no need to pay any special attention to the $y$ derivatives. They can be computed as
regular finite differences on the raw data in these frames. Thus at $o$ and in the Riemann
normal frame ay $o$
\begin{equation}
   \nabla v
      = \left(v^x \frac{1}{L_x}\frac{dL_x}{ds}\right)\partial_x
      + \frac{\partial v^y}{\partial y}\partial_y
\end{equation}

Suppose now that the vector field $v$ is the gradient of some scalar function $\phi$.
Then $v=\nabla\phi$ and
\begin{equation}
\begin{aligned}
   \nabla^2 \phi
      &= \frac{\partial v^x}{\partial x} + \frac{\partial v^y}{\partial y}\\
      &= v^y \frac{1}{L_x}\frac{dL_x}{ds} + \frac{\partial^2 \phi}{\partial y^2}\\
      &= \frac{1}{L_x}\frac{dL_x}{ds} \frac{\partial\phi}{\partial y}
       + \frac{\partial^2 \phi}{\partial y^2}
\end{aligned}
\end{equation}
and since the $y$ coordinate measures proper distance along the lattice it follows that
$d/ds = d/dy$ and thus
\begin{equation}
   \nabla^2 \phi
   = \frac{1}{L_x}\frac{dL_x}{ds} \frac{\partial\phi}{\partial s}
   + \frac{\partial^2 \phi}{\partial s^2}
\end{equation}
which agrees with (\ref{eqn:flatLaplacian}) for the particular case where $\phi=R$.

% ============================================================================================
\clearpage

\clearpage

% ============================================================================================
% \bibliographystyle{brewin}  % use this to create the bibliography
% \bibliography{brewin}

\providecommand{\href}[2]{#2}\begingroup\raggedright\endgroup
           % use this when submitting the paper to the journal
                              % or when all references have been included

\end{document}